\documentclass{amsart}

\usepackage{amssymb,amsfonts, latexsym,amsmath,hhline,array,longtable}
\usepackage{pdfsync,color, comment,colortbl, mathrsfs,stmaryrd,cite,graphicx}

\usepackage{ifpdf}
\ifpdf \usepackage[colorlinks=true, citecolor=blue, linkcolor=blue, urlcolor=blue]{hyperref} \fi

\usepackage{pdfsync}
\usepackage{xcolor} 
\usepackage{comment,colortbl}

\newcommand{\cal}{\mathcal}
\def\inp#1{\textcolor{blue}{\bf #1}}

\newtheorem{formula}{}[section]
\newtheorem{definition}[formula]{Definition}
\newtheorem{corollary}[formula]{Corollary}
\newtheorem{remark}[formula]{Remark}
\newtheorem{lemma}[formula]{Lemma}
\newtheorem{theorem}[formula]{Theorem}
\newtheorem{proposition}[formula]{Proposition}

\def\thrm{\begin{theorem}}
\def\thrml#1{\begin{theorem}\label{#1}}
\def\ethrm{\end{theorem}}
\def\rmrk{\begin{remark}}
\def\prpstn{\begin{proposition}}
\def\prpstnl#1{\begin{proposition}\label{#1}}
\def\eprpstn{\end{proposition}}
\def\rmrkl#1{\begin{remark}\label{#1}}
\def\ermrk{\end{remark}}
\def\dfntn{\begin{definition}}
\def\dfntnl#1{\begin{definition}\label{#1}}
\def\edfntn{\end{definition}}
\def\nmrt{\begin{enumerate}}
\def\enmrt{\end{enumerate}}
\def\tm#1{\item[{\rm (#1)}]}
\def\qtnl#1{\begin{equation}\label{#1}}
\def\eqtn{\end{equation}}
\def\lmm{\begin{lemma}}
\def\lmml#1{\begin{lemma}\label{#1}}
\def\elmm{\end{lemma}}
\def\crllr{\begin{corollary}}
\def\crllrl#1{\begin{corollary}\label{#1}}
\def\ecrllr{\end{corollary}}
\def\css{\begin{cases}}
\def\ecss{\end{cases}}
\def\prf{\begin{proof}}
\def\eprf{\end{proof}}

\def\cS{{\cal S}}
\def\cT{{\cal T}}

\def\mZ{{\mathbb Z}}

\def\fA{{\frak A}}

\def\fX{{\frak X}}
\def\fY{{\frak Y}}

\DeclareMathOperator{\aut}{Aut}

\DeclareMathOperator{\cay}{Cay}

\DeclareMathOperator{\cyc}{Cyc}
\DeclareMathOperator{\diag}{Diag}

\DeclareMathOperator{\hol}{Hol}
\DeclareMathOperator{\id}{id}
\DeclareMathOperator{\im}{im}
\DeclareMathOperator{\Inn}{Inn}

\DeclareMathOperator{\iso}{Iso}

\DeclareMathOperator{\orb}{Orb}

\DeclareMathOperator{\pr}{pr}

\DeclareMathOperator{\rk}{rk}

\DeclareMathOperator{\Span}{Span}

\DeclareMathOperator{\sym}{Sym}
\DeclareMathOperator{\WL}{WL}
\DeclareMathOperator{\wldim}{\dim_{\scriptscriptstyle \WL}}

\def\eprff{\hfill$\square$}

\def\grp#1{\langle {#1}\rangle}

\def\qaq{\quad\text{and}\quad}

\def\ov{\overline}
\def\phmb#1{{\phantom{x}\hspace{-2mm}^{#1}}}
\def\und#1{{\underline{#1}}}
\def\wh{\widehat}
\def\wt{\widetilde}

\begin{document}

\title{On multidimensional Schur rings of finite groups}
\author{Gang Chen}
\address{School of Mathematics and Statistics, Central China Normal University, Wuhan 430079, China}
\email{chengangmath@mail.ccnu.edu.cn}
\author{Qing Ren}
\address{School of Mathematics and Statistics, Central China Normal University, Wuhan 430079, China}
\email{renqing@mails.ccnu.edu.cn}
\author{Ilia Ponomarenko}
\address{Steklov Institute of Mathematics at St. Petersburg, Russia}
\email{inp@pdmi.ras.ru}
\thanks{}
\date{}

\begin{abstract}
For any finite group $G$ and a positive integer $m$, we define and study a Schur ring over the direct power $G^m$, which gives an algebraic interpretation of the partition of~$G^m$ obtained by the $m$-dimensional Weisfeiler-Leman algorithm. It is proved that this ring determines the group $G$ up to isomorphism if $m\ge 3$, and  approaches the Schur ring associated with the group $\aut(G)$ acting on $G^m$ naturally if  $m$ increases. It turns out that the problem of finding this limit ring is polynomial-time equivalent to the group isomorphism problem.
\end{abstract} 

\maketitle

\section{Introduction}

One way to think about the isomorphism problem of finite groups is to look for ``natural'' invariants determining a given group up to isomorphism. Examples of such invariants are group determinants~\cite{Formanek1991}, linear invariants of group rings~\cite{Roitman1981}, Cayley graphs~\cite{Baginski2021}, and so on. However, all these  invariants are hard to compute in the sense that none of the known algorithms calculates them in  time polynomial in the order of the group. On the other hand, there are many easily computable invariants that determine, up to isomorphism, not all groups but only those belonging  to a certain class. For example, abelian groups and all simple groups are uniquely determined by the multiset of orders of all elements, see, e.g., ~\cite{Grechkoseeva2021}. 

It seems quite promising to look for invariants of a  finite group $G$ among easily computable invariants for the direct powers $G^m=G\times\cdots\times G$ ($m$ times), $m=1,2,\ldots$, in the hope that for some (not too large) $m$ they determine the group~$G$ up to isomorphism. Such an approach has been introduced and studied in two recent papers \cite{Brachter2021,Brachter2020}. The idea is to define a natural canonical coloring of the elements of~$G^m$, refine the coloring with the help of the $m$-dimensional Weisfeiler-Leman algorithm ($m$-dim $\WL$), and study the invariant of~$G$ formed by numerical parameters associated with the resulting coloring. These invariants can be computed in time $n^{O(m)}$, where $n=|G|$, and determine $G$ up to isomorphism for $m=O(d)$, where here and below $d=d(G)$ denotes the minimal cardinality of a generating set for $G$.

The initial motive for writing this paper is to find an algebraic interpretation of the $m$-dim $\WL$ invariants in terms of S-rings (Schur rings). The theory of S-rings was initiated by I.~Schur (1931) and developed by H.~Wielandt, see~\cite[Chap.~IV]{Wielandt1964}; since then the S-rings are widely used in group theory and algebraic combinatorics. Recall that  a subring $\fA$  of the group ring $\mZ G$ is called an {\it S-ring} over $G$ if there exists a (uniquely determined) partition $\cS=\cS(\fA)$ of $G$ containing the identity element $1_G$ as a class, closed under taking inverse, and such that 
\qtnl{300322a}
\fA=\Span_\mZ\{\und{X}:\ X\in \cS\},
\eqtn
where $\und{X}$ denotes the sum of the elements of $X$ in $\mZ G$; the linear base of~$\fA$ consisting of the elements $\und{X}$, $X\in \cS$, is said to be \emph{standard}. In the two extreme cases, when the partition $\cS$ is discrete or consists of at most two classes ($\{1_G\}$ and its complement), the S-ring~$\fA$ is the group ring~$\mZ G$ or the \emph{trivial} ring~$\cT(G)$, respectively.

A natural example of an S-ring over $G$ (suggested by I.~Schur) is given by any permutation group $K$ containing a regular subgroup isomorphic to~$G$. In this case, the set on which $K$ acts can be identified with $G$ so that $K\le \sym(G)$ and the partition of $G$ into the orbits of the stabilizer of~$1_G$ in~$K$ is of the form $\cS(\fA)$ for some S-ring $\fA$ over~$G$. Every S-ring associated with a permutation group in this way is said to be \emph{schurian}. An obvious example of schurian S-ring is obtained for  the natural subdirect product $K=G\rtimes M$ with $M\le\aut(G)$: here,  $\cS(\fA)=\orb(M)$, and $\fA$ is denoted by $\cyc(M,G)$.  Not every S-ring $\fA$ over $G$ is schurian, but in any case, there exists a unique maximal group $K\le\sym(G)$ such that the S-ring associated with $K$ contains $\fA$ as a subring; this $K$ is called the automorphism group of~$\fA$ and denoted by $\aut(\fA)$. 

Let $G$ be a group and  $m\ge 1$ an integer.  The $m$th tensor power $\cT_m(G)$ of the trivial S-ring $\cT(G)$ is an S-ring over the group~$G^m$. If $m\ne 1$, then the ``diagonal'' element $\und{\diag(G^m)}\in\mZ G^m$ does not belong to~$\cT_m(G)$. The key ingredient of our technique is  a uniquely determined extension of $\cT_m(G)$ by this element,
$$
\fA_m(G)=\cT_m(G)[\und{\diag(G^m)}],
$$
i.e., the smallest S-ring over $G^m$, that contains both $\cT_m(G)$ and  the ``diagonal'' element.\footnote{In the category of S-rings, $\fA_m(G)$ plays the same role as the $m$-extension of a coherent configuration in the category of coherent configurations, see~\cite{EvdP1999c}. } This  S-ring is trivial for $m=1$,  and is a special case  (up to language) of the association  scheme of rank~$5$ defined  in~\cite[Subsection~2.1.5]{Heinze2009} if $m=2$ and $|G|\ge 3$. 

For $m\ge 3$, the structure of the S-ring $\fA_m(G)$ is quite complicated (though the partition $\cS(\fA_m(G))$ can easily be computed by the $2$-dim $\WL$ applying to an appropriate coloring of $G^m\times G^m$). On the other hand,  from the results of~\cite{Bailey2021} (see also Corollary~\ref{220522w}), it is not difficult to deduce that if $m\ge 2$, then
\qtnl{290322a}
\aut(\fA_m(G))=\hol_m(G),
\eqtn
where $\hol_m(G)=G^m\rtimes\aut(G)$ is the permutation group on $G^m$, induced by right multiplications of~$G^m$ and the componentwise action of $\aut(G)\le \sym(G^m)$. Moreover, we prove in Section~\ref{220522s} that the group $G$ is uniquely determined by the S-ring~$\fA_m(G)$ for $m\ge 3$.

\thrml{220522t}
For $m\ge 3$ and any two groups $G$ and $G'$, the following statements are equivalent:
\nmrt
\tm{1} the groups $G$ and $G'$ are isomorphic,
\tm{2} the S-rings $\fA_m(G)$ and $\fA_m(G')$ are combinatorially isomorphic,
\tm{3} the groups $\hol_m(G)$ and $\hol_m(G')$ are permutation isomorphic.
\enmrt
\ethrm

The equivalences in Theorem~\ref{220522t} do not hold for $m=1$ and all~$G$ of order at least~$4$. It would be interesting to characterize all $G$ for which the equivalences hold for $m=2$.

For fixed $m$, the S-ring $\fA_m(G)$ provides an easily computable invariant of the group~$G$, namely, the tensor of structure constants with respect to the standard linear base. In general, this tensor is not a full invariant, because it determines the S-ring $\fA_m(G)$ only up to algebraic (rather than combinatorial) isomorphisms  (for the exact definitions, see Section~\ref{300322c}). However, if the S-ring $\fA_m(G)$ was schurian for a fixed~$m$ and all~$G$, then the invariant would be sufficient to test isomorphism of any two groups efficiently (see Theorem~\ref{220122e} below).

Our second result shows that the sequence $\fA_1(G)$, $\fA_2(G)$, $\ldots$ is stabilized in the sense that for a sufficiently large indices,  the projections of all these S-rings to a fixed power of~$G$ become schurian. To be more precise, set $\cyc_m(G)=\cyc(\aut(G),G^m)$, and denote by $\fA_{m+k}(G)_{G^m}$, $k\ge 0$,  the image of the S-ring~$\fA_{m+k}(G)$  with respect to the natural projection of $G^{m+k}$ to the first $m$ coordinates (note that the image is an S-ring over~$G^m$).

\thrml{271221a0}
Let $G$ be a group and $m,k$ positive integers. Then  
\qtnl{300322g}
\fA_m(G)\le \fA_{m+k}(G)_{G^m}\le\cyc_m(G).
\eqtn
Moreover, the second inclusion becomes equality if $k\ge \max\{2,d\}$, where $d=d(G)$. In particular,  there exists a positive integer $t\le \max\{2,d\}$ such that
$$
\cT(G)=\ov\fA_1(G)\le \ov\fA_2(G)\le \cdots \le \ov\fA_t(G)=\ov\fA_{t+1}(G)=\cdots=\cyc_1(G),
$$
where for all~$i$, we set $\ov\fA_i(G)=\fA_i(G)_G$.
\ethrm

The S-ring $\fA_m(G)$ can naturally be interpreted in terms of  of the canonical partition $\WL_m(G)$ of $G^m$,  constructed  by the  $m$-dimensional Weisfeiler-Leman algorithm~$\WL_I$ for groups, see \cite{Brachter2020}. Namely, by Theorem~\ref{010422a} the projection of $\WL_{3m}(G)$ to $G^m$, denoted below by $\WL_{3m}(G)_{G^m}$, forms an S-ring over~$G^m$. This enables us to compare the two partitions~$\WL_m(G)$ and $\cS(\fA_m(G))$.

\thrml{230122a}
Let $G$ be a group. Then for any positive integer $m$, 
$$
\WL_{3m}(G)_{G^m}\ge \cS(\fA_m(G))\qaq\cS(\fA_{m+1}(G))_{G^m}\ge \WL_m(G).
$$
Moreover, there is an integer $s(G)\ge 3$ such that if $m\ge s(G)$, then every algebraic isomorphism from the S-ring  $\fA_m(G)$ is induced by an isomorphism, and also
\qtnl{090123a}
\wldim(G)\le 3s(G)\qaq s(G)\le \wldim(G)+1,
\eqtn
where $\wldim(G)$ is the $\WL$-dimension  of~$G$.
\ethrm

One of the primary motivations for this paper is to understand more about the group isomorphism problem, namely, test efficiently whether two given groups are isomorphic. As in the case of the graph isomorphism problem, it is quite natural to consider \emph{colored} groups: the elements of the input groups are assumed to be colored and the isomorphisms are assumed to be  color preserving, see \cite{Brachter2021}. In this paradigm, the definition of the S-ring $\fA_m(G)$ is changed in an obvious way. Namely, if $X_1,\ldots,X_c$ are colored classes of $G$, then the colored version of $\fA_m(G)$ is defined to be the extension of $\cT_m(G)$ by the elements 
$\und{\diag(X_1^m)},\ldots,\und{\diag(X_c^m)}\in\mZ G$.

The following theorem (and its proof) is very similar to the corresponding theorem (and its proof) for graphs, see~\cite{Mat1979}. However, in contrast to the case of graphs, it is not clear whether the isomorphism problems for groups and colored groups are polynomial-time equivalent. 

\thrml{220122e}
Testing isomorphism of colored groups is polynomial-time equivalent to the problem of finding the S-ring $\cyc_1(G)$ for a given colored group~$G$.
\ethrm

The paper is organized as follows. Section~\ref{300322c} provides a necessary background of S-ring theory. In Sections~\ref{020422j} and~\ref{020422w}, we describe the Weisfeiler-Leman method for groups and study the basic properties of  the S-ring $\fA_m(G)$, respectively. The proofs of the main results are given in Sections~\ref{220522s}--\ref{041122a}.

\section{Schur rings}\label{300322c}

\subsection{S-rings} Let $G$ be a finite group. A $\mZ$-submodule~$\fA$ of the group ring~$\mZ G$ is  called a {\it Schur ring} ({\it S-ring}, for short) over~$G$ if there exists a partition $\cS=\cS(\fA)$ of~$G$ such that
\nmrt
\tm{S1} $\{1_G\}\in\cS$,
\tm{S2} $X^{-1}\in\cS$ for all $X\in\cS$,
\tm{S3} $\und{X}\,\und{Y}=\sum_{Z\in\cS}c_{X,Y}^Z\und{Z}$ for all $X,Y\in \cS$ and some  integers $c_{X,Y}^Z$.
\enmrt
The elements of $\cS$ and the number $\rk(\fA)=|\cS|$ are called, respectively, the {\it basic sets} and {\it rank} of~$\fA$. The basic set containing $x\in G$ is denoted by $[x]$. The nonnegative integer $c_{X,Y}^Z$ is equal to the number of representations $z=xy$ with $x\in X$ and $y\in Y$ for a fixed $z\in Z$.

Any union of basic sets is called an \emph{$\fA$-set}. The set of all of them is closed with respect to taking inverse,  product, and standard set-theoretical operations. An $\fA$-set which is a subgroup of~$G$ is called an {\it $\fA$-group}. For example, if $X$ is an $\fA$-set, then  the group $\grp{X}$ generated by~$X$ is an $\fA$-group.  For any basic set $X$ and $\fA$-group~$H$, we put
\qtnl{180522h}
n(X,H)=\sum_Yc_{Y,X}^X,
\eqtn
where $Y$ runs over the basic sets contained in~$H$. It is not hard to prove  that $n(X,H)=|X\cap Hx|$ for every $x\in X$. 

Let $H$ be a normal $\fA$-subgroup of $G$, and $\pi:G\to G/H$ the natural epimorphism. Then the set $\cS_{G/H}=\{\pi(X):\ X\in\cS\}$ forms a partition of $G/H$. Moreover, the elements $\und{\pi(X)}$ span an S-ring over $G/H$; it is denoted by $\fA_{G/H}$ and is called the \emph{quotient} of $\fA$ modulo $H$.

The partial order  $\le$ on the  S-rings over $G$ is induced by inclusion. Thus, $\fA\le\fA'$ if and only if any basic set of $\fA$ is a union of some basic sets of $\fA'$; in this case we say that $\fA'$ is an \emph{extension} of~$\fA$. The least and greatest elements with respect to~$\le$ are, respectively,  the \emph{trivial}  S-ring $\cT(G)$ spanned by~$\und{1_G}$ and~$\und{G}$, and the group ring~$\mZ G$. For any $X\subseteq G$, the extension $\fA[\und{X}]$ of~$\fA$ by~$\und{X}$ is the smallest S-ring over~$G$, that contains~$\fA$ and~$\und{X}$.

\subsection{Isomorphisms}\label{080522f}
For any $X\in\cS$, one can define a Cayley graph $\cay(G,X)$ with vertex set $G$ so that the vertices $x$ and $y$ are adjacent if and only if $xy^{-1}\in X$. The automorphism group of this graph contains a regular subgroup of the symmetric group $\sym(G)$, induced by right multiplications of~$G$.

Let $\fA$ be an S-ring over $G$ and $\fA'$ an S-ring over $G'$. A bijection $f:G\to G'$ is called a {\it (combinatorial) isomorphism} from $\fA$ to $\fA'$ if for each $X\in\cS(\fA)$ there is $X'\in\cS(\fA')$ such that
$$
f\in\iso(\cay(G,X), \cay(G',X')),
$$ 
or, equivalently, $(Xy)^f=X'f(y)$ for all $y\in G$. In particular, the set $\iso(\fA,\fA')$ of all isomorphisms from $\fA$ to $\fA'$ contains all group isomorphisms $f:G\to G'$ for which $\cS(\fA)^f=\cS(\fA')$. Any such~$f$ is a \emph{normalized} isomorphism, i.e., $f(1_{G^{}})=1_{G'}$. Note that if the S-rings~$\fA$ and~$\fA'$ are isomorphic if and only if there is a normalized isomorphism from~$\fA$ to~$\fA'$.

The group $\iso(\fA)=\iso(\fA,\fA)$ of all  isomorphisms from $\fA$ to itself has a normal subgroup equal to the intersection of the groups $\aut(\cay(G,X))$, $X\in\cS$. It is called the \emph{automorphism group} of~$\fA$ and denoted by~ $\aut(\fA)$. In particular,  
\qtnl{311222a}
\aut(\fA)\le\aut(\cay(G,X))
\eqtn
for each~$X$. Thus, $f$ is an automorphism of~$\fA$ if and only if $(Xy)^f=Xy^f$ for all $X\in\cS$ and all $y\in G$. Note that $f\in\aut(G)$ belongs to $\aut(\fA)$ if and only if $f$ leaves every basic set of~$\fA$ fixed. 

\lmml{140514a}
Let $\fA$ and $\fA'$ be S-rings over $G$ and $G'$, respectively, and let $X\subseteq G$ and $X'\subseteq G'$.  Suppose that $f\in\iso(\fA,\fA')$ is a normalized isomorphism such that  
$$
f(Xy)=X'f(y)
$$
for all $y\in G$. Then $f\in\iso(\fA[\und{X}],\fA'[\und{X'}])$. In particular, if $G=G'$, $\fA=\fA'$, and $X=X'$, then $f\in\aut(\fA[\und{X}])$.
\elmm
\prf Follows from \cite[Lemma~2.3]{Evdokimov2015}.
\eprf

Let $m\ge 1$ be an integer, $A\le\aut(G)$, and $\cS$ the partition of~$G^m$ into the orbits of the componentwise action of~$A$. Then the module $\fA$ defined by formula~\eqref{300322a} is an S-ring over~$G^m$. It is denoted by $\cyc_m(A,G)$. This is a particular example of an S-ring associated with subgroup of $\sym(G)$, containing the permutations induced by right multiplications of~$G$. 

\subsection{Algebraic isomorphisms}
In the notation of the previous subsection, a ring isomorphism $\varphi:\fA\to\fA'$ is called an \emph{algebraic isomorphism} if for any $X\in\cS$ there exists $X'\in\cS'$ such that $\varphi(\und{X})=\und{X'}$. From the definition, it follows  that the mapping $X\mapsto X'$ is a bijection from $\cS$ onto $\cS'$. This bijection is naturally extended to a bijection between the $\fA$- and $\fA'$-sets, that takes the $\fA$-groups to $\fA'$-groups;  the image $X'$ of an $\fA$-set $X$ is also by $\varphi(X)$. The equalities 
$$
c_{X,Y}^Z=c_{X',Y'}^{Z'}\qaq n(X,H)=n(X',H')
$$
hold for all basic sets $X,Y,Z$ and for all $\fA$-groups $H$. If $H$ is a normal $\fA$-subgroup of~$G$,  $\pi:G\to G/H$ the natural epimorphism, and $H'$ a normal subgroup of $G'=\varphi(G)$, then the mapping $\pi(X)\mapsto \pi'(X')$, $X\in\cS$, induces an algebraic  isomorphism $\varphi_{G/H}:\fA^{}_{G/H}\to\fA'_{G'/H'}$ such that
\qtnl{180522a}
\varphi_{G/H}(\pi(X))=\pi'(\varphi(X)),
\eqtn
where $\pi':G'\to G'/H'$ is the natural epimorphism.

\lmml{100522b}
Let $\varphi:\fA\to\fA'$ be an algebraic isomorphism. Then for any S-rings $\wt\fA\ge \fA$ and $\wt\fA'\ge\fA'$ there is at most one algebraic isomorphism from $\wt\fA$ to $\wt\fA'$  extending $\varphi$.
\elmm
\prf Follows from \cite[Lemma~2.1]{Evdokimov2015}.\eprf

Every normalized isomorphism $f:\fA\to\fA'$ defines a bijection $\varphi:\cS\to\cS'$, $X\mapsto X^f$, which is  an algebraic isomorphism; we say that $\varphi$ is induced by~$f$. For example, the trivial algebraic isomorphism $\varphi=\id$ is induced by every normalized automorphism of~$\fA$. However, not every algebraic isomorphism is induced by an isomorphism.

\subsection{Tensor product} Let $\fA$ and $\fA'$ be  S-rings over groups $G$ and $G'$, respectively. 
Then the Cartesian products $X\times X'$ with $X\in\cS$ and $X'\in\cS'$, form a partition $\cS\times\cS'$ of the direct product $G\times G'$. Moreover, there is a uniquely determined S-ring $\fA\otimes\fA'$ over $G\times G'$ such that 
$$
\cS(\fA\otimes\fA')=\cS\times\cS'.
$$
It is called the {\it tensor product} of $\fA$ and $\fA'$. The following lemma was proved in \cite[Lemma~2.2]{Evdokimov2013}.

\lmml{170422a}
Let $G$ and $G'$ be groups  and $\pi:G\times G'\to G$ and $\pi':G\times G'\to G'$ the natural projections. Let $\fA$ be an S-ring over $G\times G'$ such that $G$ and $G'$ are $\fA$-groups. Then $\pi(X),\pi'(X)\in\cS(\fA)$ for all $X\in\cS(\fA)$.  In particular, $\fA\ge\fA_{G^{}}\otimes \fA_{G'}$. 
\elmm

The tensor product of $m$ copies of the trivial S-ring $\cT_G$ is denoted by $\cT_m(G)$. The basic sets $X$ of this S-ring are in one-to-one correspondence with the sets $I\subseteq\{1,\ldots,m\}$; namely, $X=\{x\in G^m:\ x_i=1_G$ if and only if $i\in I\}$.

\section{The Weisfeiler-Leman method for groups}\label{020422j}

The key point in the Weisfeiler-Leman ($\WL$) method for groups is the multidimensional $\WL$ algorithm. For the purpose of the present paper, it is not necessarily to describe this algorithm in detail, and it suffices to know the structure of the resulted multidimensional coherent configuration introduced in~\cite{Babai2019}, see also~\cite{AndresHelfgott2017}. 

\subsection{Multidimensional coherent configurations}\label{140422a}
Let $\Omega$ be a finite set  and $m$ a positive integer. Let us fix some notation. For a tuple $x\in\Omega^m$, denote by $\rho(x)$ the equivalence relation on $M=\{1,\ldots,m\}$ such that $(i,j)\in\rho(x)$ if and only if  $x_i=x_j$.  The class of  a partition $\fX$ containing $x$,  is denoted by~$[x]$.
For a mapping $\sigma:M\to M$, we put  $x^\sigma=(x_{1^\sigma},\ldots,x_{m^\sigma})$.

\dfntn
A partition $\fX$ of  $\Omega^m$ is called an $m$-ary rainbow on $\Omega$  if the following conditions are satisfied for all $X\in \fX$:
\nmrt
\tm{C1} $\rho(x)$ does not depend on $x\in X$,
\tm{C2} $X^\sigma\in \fX$  for all mapping $\sigma:M\to M$.
\enmrt
\edfntn

The $m$-ary rainbows on $\Omega$ are (up to language) the $m$-ary configurations in~\cite{AndresHelfgott2017}: the  difference is that we do not use a coloring of~$\Omega^m$ to define the classes. Condition~(C2) implies that the coordinatewise action of~$\sym(m)$ on~$\Omega^m$ induces an action on~$\fX$; in particular, 
$$
[x]^\sigma\in \fX\quad\text{for all } x\in\Omega^m,\ \sigma\in\sym(m).
$$
The $2$-ary rainbows are ordinary rainbows in the sense of~\cite{CP2019}, but the converse statement is not necessarily true.

Let $\fX'$ be an $m$-ary rainbow on $\Omega'$. A bijection $f:\Omega\to\Omega'$ is called an \emph{isomorphism} from $\fX$ to~$\fX'$ if the induced bijection 
$$
f_m:\Omega^m\to{\Omega'}\phmb{m},\ (x_1,\ldots,x_m)\mapsto (x_1^f,\ldots,x_m^f),
$$ 
takes the classes of $\fX$ to those of $\fX'$. Clearly, $\rho(x)=\rho(x^{f_m})$ for all $x\in\Omega^m$, and $X^{f_m\sigma}=X^{\sigma f_m}$ for all  mappings $\sigma:M\to M$.

The $m$-ary rainbows on $\Omega$ are partially ordered in accordance with the partial order of partitions. Namely,  $\fX\le \fY$ if every class of $\fX$ is a union of some classes of~$\fY$, or equivalently, if $\fX^\cup\subseteq \fY^\cup$, where $\fX^\cup$ (respectively, $\fY^\cup$) is the set of all unions of classes of~$\fX$ (respectively,~$\fY$). The largest $m$-ary rainbow is the \emph{discrete} one in which every class is a singleton; the smallest $m$-ary rainbow consists of the orbits of the symmetric group $\sym(\Omega)$ in its componentwise action on~$\Omega^m$. Given $m$-ary rainbows $\fX$ and $\fY$ on the same set, there is a unique $m$-ary rainbow $\fX\cap\fY$ such that $(\fX\cap\fY)^\cup=\fX^\cup\cap\fY^\cup$; it is called the \emph{intersection} of~$\fX$ and~$\fY$.

Let $\fX$ be an $m$-ary rainbow on $\Omega$. For any $x\in\Omega^m$ and any $X_1,\ldots,X_m\in \fX$, denote by $n(x;X_1,\ldots,X_m)$ the number of all $\alpha\in\Omega$ such that $x_{i\leftarrow \alpha} \in X_i$ for all~$i\in M$, where
\qtnl{030422q}
x_{i\leftarrow \alpha}=(x_1,\ldots,x_{i-1},\alpha,x_{i+1},\ldots,x_m).
\eqtn
These numbers define an equivalence relation $\sim$ on~$\Omega^m$, such that $x\sim y$ if and only if for all $X_1,\ldots,X_m\in\fX$,
$$
[x]=[y]\qaq n(x;X_1,\ldots,X_m)= n(y;X_1,\ldots,X_m).
$$ 
The classes of $\sim$ form a partition of~$\Omega^m$, called the \emph{$\WL$-refinement} of $\fX$ and denoted by $\WL_{m,1}(\fX)$. It is not hard to verify that the $\WL$-refinement  takes $\fX$ to a partition satisfying~(C1). 

For a given $\fX$, the $m$-dim Weisfeiler-Leman algorithm $\WL_m$ constructs  the following partitions step by step: 
\qtnl{130522q}
\fX=\fX_0< \fX_1<\ldots <\fX_k=\fX_{k+1}=\WL_m(\fX)
\eqtn
for some $k\le |\Omega|^m$, where $\fX_{i+1}=\WL_{m,1}(\fX_i)$ for $i=0,1,\ldots,k$. The resulted partition  $\WL_m(\fX)$ is the smallest $m$-ary coherent configuration on $\Omega$ in the sense of the definition below that is larger than or equal to~$\fX$. Note that the mapping $\fX\mapsto \WL_m(\fX)$ is easily computable and defines a closure operator on the $m$-ary rainbows. In particular, it is monotone, i.e.,  $\fX\le \fX'$ implies $\WL_m(\fX)\le  \WL_m(\fX')$.

\dfntn
An $m$-ary coherent configuration is  an $m$-ary rainbow $\fX$  satisfying the additional  condition
\nmrt
\tm{C3} for any $X_0,X_1,\ldots,X_m\in \fX$, the number $n_{X_1,\ldots,X_m}^{X_0}=n(x_0;X_1,\ldots,X_m)$ 
does not depend on $x_0\in X_0$.
\enmrt
\edfntn

The unary coherent configurations on~$\Omega$ are just the partitions of~$\Omega$, whereas the $2$-ary coherent configurations are ordinary coherent configurations, see, e.g., \cite{CP2019}. A natural example of an $m$-ary coherent configuration is given by the set of orbits of the natural action on~$\Omega^m$ of a permutation group on~$\Omega$.

Following~\cite[Formula~(9)]{Ponomarenko2022a}, an \emph{algebraic isomorphism} of $m$-ary coherent configurations $\fX$ and $\fX'$ is a bijection $\varphi:\fX\to\fX'$ such that for all $X,X_0,\ldots,X_m\in \fX$ and $\sigma:M\to M$,
\qtnl{110522a}
\varphi(X^\sigma)=\varphi(X)^\sigma\qaq n_{X_1,\ldots,X_m}^{X_0}=n_{\varphi(X_1),\ldots,\varphi(X_m)}^{\varphi(X_0)}.
\eqtn
Every isomorphism $f$ from $\fX$ to $\fX'$ induces an algebraic isomorphism $\varphi:\fX\to\fX'$ such that $X^\varphi=X^{f_m}$ for all~$X$. Furthermore, if   $X\in \fX$ and $i,j\in M$, then $(i,j)\in\rho(X)$ if and only if $X^\sigma=X$, where $\sigma:M\to M$ is identical on $M\setminus\{j\}$ and takes $j$ to~$i$. It follows that
\qtnl{110622x}
\rho(\varphi(X))=\rho(X).
\eqtn
We extend $\varphi$ to a bijection $\fX^\cup\to (\fX')^\cup$ in a natural way. Then $X\subseteq Y$ implies $\varphi(X)\subseteq \varphi(Y)$ for all $X,Y\in\fX^\cup$.

\lmml{130522v}
Let $\varphi$ be an algebraic isomorphism from an $m$-ary coherent configuration $\fX$ to another  $m$-ary coherent configuration. Then for every rainbow $\fY\le \fX$, 
$$
\varphi(\WL_m(\fY))=\WL_m(\varphi(\fY)).
$$
\elmm
\prf 
Using induction on $i=0,\ldots,k$ in formula~\eqref{130522q} for $\fX=\fY$, it suffices to prove that $\varphi(\WL_{m,1}(\fY_i))=\WL_{m,1}(\varphi(\fY_i))$. In other words, we need to verify that if $Y\in\fY_i$ and $Y'\in\WL_{m,1}(\fY_i)$, then 
$$
Y'\subseteq Y\quad\Rightarrow\quad \varphi(Y')\subseteq\varphi(Y).
$$
But this easily follows from formulas~\eqref{110522a} applied to the classes of~$\fX$, that are contained in the classes of~$Y'$.\eprf

\subsection{Projections}
Let $K=\{i_1,\ldots,i_k\}$, where $1\le i_1<i_2<\ldots<i_k\le m$. The projection of $\Omega^m$ with respect to a set $K\subseteq M$ is  the mapping 
$$
\pr_K:\Omega^m\to \Omega^k,\ (\alpha_1,\ldots,\alpha_m)\mapsto (\alpha_{i_1},\ldots,\alpha_{i_k}).
$$
When $K=\{1,\ldots,a\}$, we abbreviate $\pr_a:=\pr_K$. 

\lmml{090122b}
Let $\fX$ be an $m$-ary coherent configuration (respectively, rainbow) on~$\Omega$, and $K\subseteq M$. Then $\pr_K(\fX)=\{\pr_K(X):\ X\in \fX\}$ is a $|K|$-ary coherent configuration (respectively, rainbow) on $\Omega$. In particular,
\qtnl{160523g}
\pr_K^{-1}(\pr_K(X))\in \fX^\cup\quad \text{ for all }X\in\fX^\cup.
\eqtn
Moreover, every algebraic isomorphism $\varphi:\fX\to\fX'$ induces an algebraic isomorphism $\varphi_K:\pr_K(\fX)\to\pr_K(\fX')$ such that
\qtnl{160523h}
\varphi_K(\pr_K(X))=\pr_K(\varphi(X))\quad \text{ for all }X\in\fX^\cup.
\eqtn
\elmm

\prf
The first statement was proved  in  \cite[Exercises 2.7, 2.11]{AndresHelfgott2017} (see also \cite[Lemma~3.3]{Ponomarenko2022a}), the second one is trivial and the third is  
\cite[Lemma~3.3]{Ponomarenko2022a}.
\eprf

Let $\fX$ be an $m$-ary rainbow on $\Omega$. For every $X\in \fX$, $K\subseteq M$, and $x\in X$, denote by $n_K(x; X)$ the number of all $y\in X$ such that $\pr_K(y)=\pr_K(x)$. When $\fX$ is an $2$-ary coherent configuration, this number is just the valency of~$X$ and does not depend on~$x\in X$. The following statement generalizes this property to arbitrary~$m$ (for some special $m$-ary coherent configurations, this was proved in~\cite[Theorem~6.1]{EvdP1999c}).

\lmml{090122d}
Let $\fX$ be an $m$-ary coherent configuration, $X\in \fX$, and $K\subseteq M$. Then the number $n_K(X)=n_K(x;X)$ does not depend on $x\in X$.
\elmm
\prf Let $x\in X$ and $A=\{y\in X:\ \pr_K(y)=\pr_K(x)\}$. First, assume that $K=\{1,\ldots,m-1\}$. When the classes $X_1,\ldots,X_{m-1}$ run over $\fX$, the nonempty  sets
$$
A(X_1,\ldots,X_{m-1})=\{y\in A:\ x_{1\leftarrow y_m}\in X_1,\ldots,x_{m-1\leftarrow y_m}\in X_{m-1}\}
$$
form a partition of~$A$. Furthermore, $|A(X_1,\ldots,X_{m-1})|=n_{X_1,\ldots,X_{m-1},X}^X$. Thus,
the number 
$$
n_K(x;X)=|A|=\sum_{X_1,\ldots,X_{m-1}}|A(X_1,\ldots,X_{m-1})|=\sum_{X_1,\ldots,X_{m-1}}n_{X_1,\ldots,X_{m-1},X}^X
$$
does not depend on~$x\in X$.

Now, without loss of generality, we may assume that $K=\{1,\ldots,k\}$ for some $1\le k< m-1$.  Let $K'=\{1,\ldots,m-1\}$.  By Lemma~\ref{090122b}, $\fX'=\pr_{K'}(\fX)$ is an $(m-1)$-ary coherent configuration. Put  $A'=\pr_{K'}(A)$,  $X'=\pr_{K'}(X)$, and $x'=\pr_{K'}(x)$. By induction, we have
$$
n_K(X')=n_K(x';X')=|A'|.
$$
Note that if $y\in A$, then $y'\in A'$. Moreover, by the statement proved in the first paragraph, there are exactly $n_{K'}(y)=n_{K'}(X)$ elements $z\in X$ for which $\pr_{K'}(z)=y'$. Thus, the number
$$
|A|=n_K(x;X)=|A'|\,n_{K'}(x)=n_K(X')\,n_{K'}(X)
$$
does not depend on~$x\in X$.\eprf

An $m$-ary rainbow $\fX$ is said to be \emph{regular} if the statement of Lemma~\ref{090122d} holds for every $X\in\fX$ and all $K\subseteq M$. Thus Lemma~\ref{090122d} states that every $m$-ary coherent configuration is regular. The converse statement is not true even for $m=2$.
 
\crllrl{261222a}
Let $\varphi$ be an algebraic isomorphism from an $m$-ary coherent configuration $\fX$, then  $|X|=|\varphi(X)|$ for all $X\in\fX^\cup$.
\ecrllr

\subsection{The Weisfeiler-Leman partition of a direct power} 
Let $G$ be a group and $m\ge 1$.  Denote by $\fX=\fX_m(G)$ the minimal $m$-ary rainbow on $G$ such that for every $X\in \fX$ and every $x,y\in X$, the equivalences
\qtnl{030422a}
x_i=x_j\ \Leftrightarrow\ y_i=y_j\qaq x_ix_j=x_k\ \Leftrightarrow\ y_iy_j=y_k
\eqtn
hold true for all $i,j,k\in M$. The first part means exactly that the equivalence relation $\rho(x)$ does not depend on $x\in X$, whereas the second one means that the same is true for the ternary relation 
$$
\mu(x)=\{(i,j,k)\in M^3:\ x_ix_j=x_k\}.
$$ 
The  partition of~$G^m$ defined by formulas~\eqref{030422a} was used in \cite{Brachter2020} as the initial coloring of the group $G^m$ for the $\WL$ algorithm for groups in version~I.  

\lmml{140422b}
Let $G$ be a group and $\fX=\fX_m(G)$. Then
\nmrt
\tm{1} $\{1_{G^m}\}\in \fX$  and $\diag(G^m)\in \fX^\cup$,
\tm{2} $\fX^{-1}=\fX$,
\tm{3} $\fX\ge\cS(\cT_m(G))$.
\enmrt
\elmm
\prf  Statement~(1) immediately follows from~\eqref{030422a}. Next, for every $x\in G^m$, we have $\rho(x)=\rho(x^{-1})$ and $(i,j,k)\in\mu(x)$ if and only if $(j,i,k)\in\mu(x^{-1})$. This proves statement~(2). Finally, statement~(3) holds by the remark at the end of Section~\ref{300322c}.
\eprf

For a  group $G$, the smallest $m$-ary coherent configuration containing $\fX_m(G)$ is denoted by $\WL_m(G)$,
$$
\WL_m(G)=\WL_m(\fX_m(G)).
$$
It coincides with the partition obtained by the algorithm $\WL_m^I$ applied to~$G$ for $m\ge 2$, see~\cite{Brachter2020}. Note that $\WL_1(G)$ is a partition in at most two classes.

Let $\fX\ge  \WL_m(G)$ and $\fX'\ge \WL_m(G')$ be $m$-ary coherent configurations. The algebraic isomorphism $\varphi:\fX\to\fX'$  is said to be \emph{genuine} if for all $X\in\fX$,
\qtnl{100522j}
\mu(X)=\mu(\varphi(X)).
\eqtn
The groups $G$ and $G'$ are said to be \emph{$\WL_m$-equivalent} if there exists a genuine algebraic isomorphism  $\varphi:\WL_m(G)\to\WL_m(G')$.  This concept corresponds to the equivalence of $G$ and $G'$ with respect to $m$-$\WL$ in version I, introduced in~\cite{Brachter2020} for $m\ge 2$. Clearly, any two isomorphic groups are $\WL_m$-equivalent for all~$m$. The \emph{$\WL$-dimension} $\wldim(G)$ of the group~$G$ is defined to be the smallest $m$ for which every group $\WL_m$-equivalent to~$G$ is isomorphic to~$G$. For more details, we refer to~\cite{Brachter2020}.

\section{The S-ring $\fA_m(G)$}\label{020422w}

Throughout the rest of the paper, $G$ is a finite group, $m$ a positive integer, $M=\{1,\ldots,m\}$, $\fX_m=\fX_m(G)$, $\cT_m=\cT_m(G)$, $D_m=\diag(G^m)$, $\fA_m=\fA_m(G)$, and $\cS_m=\cS(\fA_m)$.

\subsection{Basic properties}  Recall that the S-ring $\fA_m$ is the extension of the $m$th tensor power~$\cT_m$ by the element $\und{D_m}$. Every $\sigma\in\sym(M)$ induces (by permuting coordinates) a normalized isomorphism of the S-ring $\cT_m$. Moreover, $(D_mx)^\sigma=D_mx^\sigma$ for all $x\in G^m$. By Lemma~\ref{140514a}, this implies that $\sigma$ induces an isomorphism of $\fA_m$ to itself. In particular, we have the following statement.

\lmml{170422c}
$X^\sigma\in\cS_m$ for all $X\in \cS_m$ and $\sigma\in\sym(M)$.
\elmm

Let $K\subseteq M$ and $K'=M\setminus K$.  Put 
$$
G_K=\{x\in G^m:\ \pr_{K'}(x)=1_{G^{K'}}\},
$$
and abbreviate $G_i=G_{\{i\}}$ and $G_{i'}=G_{M\setminus\{i\}}$ for $i\in M$. It is easily seen that~$G_K$ is an $\cT_m$-group and hence an $\fA_m$-group (because $\fA_m\ge \cT_m$). Thus, $G^m$ is the direct product  of the~$\fA_m$-groups~$G_{K^{}}$ and~$G_{K'}$. Under the factorization of $\fA_m$ modulo~$G_{K'}$, the diagonal group $D_m$ goes to $D_k$, where $k=|K|$.  By  Lemma~\ref{170422a},  the full preimage
$$
D_K=D_K(G) =\{x\in G^m:\ K^2\subseteq \rho(x)\} 
$$  
of the group~$D_k$ is  an $\fA_m$-group. Note that $D_m=D_K$ for $K=M$. Thus, we proved the following lemma.

\lmml{170422d}
$G_K$ and  $D_K$  are $\fA_m$-groups for all $K\subseteq M$.
\elmm

In order to compare the partitions $\cS_m$ and $\WL_m(G)$, we define the set
$$
X_{i,j,k}(G)=\{x\in G^m:\ (i,j,k)\in\mu(x)\},\qquad i,j,k\in M.
$$

\crllrl{250422d}
If  $m\ge 3$ and $i,j,k\in M$, then $X_{i,j,k}(G)$ is an $\fA_m$-set.
\ecrllr
\prf
If $i\ne j\ne k\ne i$, then $X_{i,j,k}(G)=(D_{\{i,k\}}\cap G_{j'})\,\cdot\,(D_{\{j,k\}}\cap G_{i'})$ and we are done by Lemma~\ref{170422d}. In the remaining case $X_{i,j,l}=(X_{i',j',k'})^\sigma$ for suitable pairwise distinct $i',j',k'\in M$ and a mapping $\sigma:M\to M$ taking $i,j,k$ to $i',j',k'$, respectively.
\eprf

Let $k\in M$ and $K=\{1,\ldots,k\}$. The group $G_K$ can naturally be identified with direct power $G^k$, and the quotient of  the  S-ring $\fA_m$ modulo~$G_{K'}$ can naturally be identified with the S-ring $\fA=(\fA_m)_{G^k}$ defined in the introduction.  Under taking this quotient, $\cT_m$ and $D_m$ go to $\cT_k$ and $D_k$, respectively. It follows that $\fA\ge \cT_k$ and $D_k$ is an $\fA$-set. Consequently, $\fA$ contains the extension of $\cT_k$ by~$D_k$, which is just $\fA_k$. Thus, \qtnl{180422d}
\fA_m(G)_{G^k}\ge \fA_k(G).
\eqtn

The following statement enables us to interpret the partition $\cS_m$ in terms  used in Section~\ref{020422j}. 

\lmml{160422w}
For every $m\ge 1$, the partition $\cS_m$ is a regular $m$-ary rainbow.
\elmm
\prf Let $X\in\cS_m$. Let $i,j\in M$ and $x\in X$ be such that $(i,j)\in\rho(x)$.  Then $x\in X\cap D_{\{i,j\}}$. 
By Lemma~\ref{170422d}, this implies that $X\subseteq D_{\{i,j\}}$. Hence, $(i,j)\in\rho(x)$ for  all $x\in X$. This proves condition~(C1).

To verify condition (C2), let $\sigma:M\to M$ be an arbitrary mapping. It can be written as a composition of a  permutation of~$M$ and  some mappings $\sigma_{i,j}:M\to M$, where $i,j\in M$, such that $j^{\sigma_{i,j}}=i$ and $k^{\sigma_{i,j}}=k$ for $k\ne j$. By Lemma~\ref{170422c}, it suffices to prove that $X^{\sigma_{i,j}}\in\cS_m$ for all~$i,j$. But this follows from the obvious equality
$$
X^{\sigma_{i,j}}=XG_i\,\cap\, D_{\{i,j\}}.
$$

Thus, $\cS_m$ is an $m$-ary rainbow. It remains to verify, that for every $K\subseteq M$ the number $n_K(x,X)$ does not depend on $x\in X$.  However,  this is true, because $n_K(x,X)= |X\cap G_{K'}x|=n(X,G_{K'})$, see formula~\eqref{180522h}. \eprf

It is not clear whether $\cS_m$ is always an $m$-ary coherent configuration for all~$m$. This is obvious for $m=1$ and we have no counterexample for $m>1$.

Recall that $\fX_m$ is the smallest  $m$-ary rainbow on $\Omega$, such that  for every $X\in\fX_m$ the  ternary relation $\mu(x)$ does not depend on $x\in X$. On the other hand, $\cS_m$ is an $m$-ary rainbow on $\Omega$ by Lemma~\ref{160422w} and satisfies this property by Corollary~\ref{250422d}. Thus the following statement holds.

\crllrl{180422b}
$\cS_m\ge \fX_m$ for all~$m\ge 3$.
\ecrllr

For $m=2$, the statement of Corollary~\ref{180422b} does not hold. Indeed, let $G_1=G\times 1_G$, $G_2=1_G\times G$, and $G_3=\diag(G\times G)$. Then
$$
\und{G_i}\cdot\und{G_j}=\und{G\times G},\qquad 1\le i\ne j\le 3.
$$
It follows that if $G_0=\{1_{G^2}\}$ and $X$ is the complement of $G_1\cup G_2\cup G_3$, then the elements $\und{G_0},\ldots,\und{G_3}$ and  $\und{X}$ span an S-ring $\fA$ over $G^2$ of rank~$5$.  It is easily seen that $\cS_2=\cS(\fA)$, and hence $|\cS_2|=5$. On the other hand, if $G$ is a cyclic $2$-group and~$x$ is a unique involution of $G$, then $\{(1_G,x)\}\in \fX_2$. Therefore, $|\fX_2|>5$ if $|G|\ge 8$. Thus, $\cS_2\not \ge \fX_2$.

\subsection{Classes of $\cS_m$} Let $w$ be an \emph{$m$-word} by which we mean a word in the alphabet consisting the $2m$ letters~$a_1^{\pm 1},\ldots,a_m^{\pm 1}$. One can define a word map  $G^m\to G$ taking an $m$-tuple $x$ to the evaluation $w(x)\in G$ replacing $a^{\pm 1}_i$ with $x^{\pm 1}_i$ for all $i\in M$ (if $w$ is an empty word, then  $w(x)=1_G$ for all~$x$). For example, if  $m=2$ and $w=a_1a_2$, then the above map defines the multiplication table of the group~$G$.  

\thrml{180422a1}
Let $0\le k\le  m-2$ and $k+1\le \ell\le m$.  Then for every $k$-word $w$ and every $X\in\cS_m$, the equality $x_\ell=w(x_1,\ldots,x_k)$ holds for all or for no $x\in X$.
\ethrm
\prf  Induction on the length $|w|$ of the word~$w$. Let $X\in\cS_m$. If $|w|=0$, then~$w$ is empty and the equality $x_\ell=w(x_1,\ldots,x_k)$ for some $x\in X$ exactly means that $x_\ell=1_G$ or equivalently, $x\in  G_{\ell'}$. Since $G_{\ell'}$ is an $\fA_m$-group, this shows that $X\subseteq G_{\ell'}$, i.e., $x_\ell=w(x_1,\ldots,x_k)$ for all $x\in X$.  

Let  $|w|\ge 1$. Then $w=w' a^{}_i$ or $w'a^{-1}_i$ for  some $k$-word $w'$  and $1\le i\le k$. We consider the first case; the second one is similar. We may assume that $\ell=m-1$ (Lemma~\ref{170422c}) and  $x_\ell=w(x_1,\ldots,x_k)$ for some $x\in X$. Put $X'=[x']$, where
\qtnl{290422a}
x'=(x_1,\ldots,x_k,\ldots,x_{m-2},w'(x_1,\ldots,x_k),w(x_1,\ldots,x_k)).
\eqtn
Since $|w'|=|w|-1$, the induction hypothesis implies that $y_{m-1}=w'(y_1,\ldots,y_k)$ for all $y\in X'$. Furthermore, 
$$
x'_{m-1}x'_i=w'(x_1,\ldots,x_k)x'_i=w(x_1,\ldots,x_k)=x'_m.
$$ 
By
Corollary~\ref{250422d}, we have $X'\subseteq X_{m-1,i,m}(G)$. Thus if $\sigma=(m-1,m)\in\sym(M)$, then every element of the set $Y=(X')^\sigma G_m$ is of the form
$$
y=(y_1,\ldots,y_{m-2},w(y_1,\ldots,y_k),y_m).
$$
On the other hand, $Y$ is an $\fA_m$-set by Lemma~\ref{170422d}, and $x\in Y$. Thus, $X\subseteq Y$ and the equality $x_{m-1}=w(x_1,\ldots,x_k)$ holds for all $x\in X$. \eprf

\crllrl{290422f}
Let $m\ge 3$ and $0\le k\le m-2$. Assume that $x\in X\in\cS_m$ is such that
\qtnl{290422e}
\grp{x_1,\ldots,x_m}=\grp{x_1,\ldots,x_k}.
\eqtn
Then given $y\in X$, the mapping $x_i\mapsto y_i$, $i=1,\ldots m$, induces an isomorphism from $\grp{x_1,\ldots,x_m}$ to $\grp{y_1,\ldots,y_m}$.
\ecrllr
\prf 
By Theorem~\ref{180422a1}, it suffices to verify that the mapping $x_i\mapsto y_i$, $i=1,\ldots k$, induces an isomorphism from $\grp{x_1,\ldots,x_k}$ to $\grp{y_1,\ldots,y_k}$. However, this is true if 
\qtnl{170822a}
w(x_1,\ldots,x_k)=w'(x_1,\ldots,x_k)\quad\Leftrightarrow\quad w(y_1,\ldots,y_k)=w'(y_1,\ldots,,y_k)
\eqtn
for all $k$-words $w$ and $w'$. To verify this equivalence, it suffices to prove the implication~$\Rightarrow$ only, because equality~\eqref {290422e} holds true for all $y\in X$ by Theorem~\ref{180422a1}. Let $X'=[x']$, where
$$
x'=(x_1,\ldots,x_{m-2}, w(x_1,\ldots,x_k),w'(x_1,\ldots,x_k)).
$$
The left-hand side of \eqref{170822a} implies that $(m-1,m)\in\rho(X')$. On the other hand, $\pr_{m-2}(y)\in\pr_{m-2}(X)=\pr_{m-2}(X')$ and hence one can find  $y'\in X'$ such that $\pr_{m-2}(y')=\pr_{m-2}(y)$. By Theorem~\ref{180422a1}, this implies that
$$
y'_{m-1}=w(y_1,\ldots,y_k),\quad y'_m=w'(y_1,\ldots,y_k).
$$
Since $(m-1,m)\in\rho(X')$, we conclude that $y'_{m-1}=y'_m$ implying the  right-hand side of \eqref{170822a}. 
\eprf

\section{Isomorphisms of the S-ring $\fA_m(G)$}\label{220522s}

\subsection{Diagonal groups} Our description of algebraic and combinatorial isomorphisms of the S-ring $\fA_m(G)$ is based on the main results on diagonal groups in~\cite{Bailey2021}. The following statement is a special case of \cite[Theorem~1.1(b)]{Bailey2021}.

\lmml{220522a}
Let $H$ be a group and $m\ge 3$ an integer. Assume that  $H$ is the direct product of any $m$ of subgroups $H_0,H_1,\ldots, H_m$. Then there is a group $T$   such that 
$$
H=T^m,\quad H_0=D_m(T),\quad H_i=T_i \text{ for all } i\in M.
$$
\elmm

For $m\ge 2$, the diagonal graph $\Gamma_D(G,m)$ is defined to be the Cayley graph $\cay(G^m,X_m)$, where $X_m=X_m(G)=G_0\cup G_1\cup\ldots\cup G_m$ with $G_0=D_m(G)$.  It was proved in~\cite[Subsection~7.2]{Bailey2021} that $\Gamma_D(G,m)$ determines $G$ up to isomorphism, 
\qtnl{220522c}
\Gamma_D(G,m)\cong\Gamma_D(G',m')\quad\Leftrightarrow\quad G\cong G'\text{ and }m=m',
\eqtn
and also that except for four small cases, the automorphism group of $\Gamma_D(G,m)$  is equal to  the semidirect product $D(G,m)=\hol_m(G)\rtimes\sym(m+1)$, see  \cite[Theorem~1.4(b)]{Bailey2021}. The statement below is an immediate  consequence of the description of the group $D(G,m)$ given in~\cite[Remark~1.3]{Bailey2021}. In what follows, we denote by $\aut_m(G)$ the permutation group induced by the componentwise action of $\aut(G)$ on $G^m$.

\lmml{220522b}
Let $G$ be a group and $m\ge 2$ an integer. Denote by $D_0(m,G)$ the subgroup of $\aut(\Gamma_D(G,m))$, that leaves the vertex $1_{G^m}$ and the each of the  sets $G_0,G_1,\ldots,G_m$ fixed.   Then $D_0(m,G)=\aut_m(G)$.
\elmm

\crllrl{220522w}
For $m\ge 2$, we have $\aut(\fA_m(G))=\hol_m(G)$.
\ecrllr
\prf Set $\fA=\fA_m(G)$ and $\ov M=M\cup\{0\}$. Recall that $G_i$ is an $\fA$-group for every $i\in \ov M$. Hence the stabilizer $\aut(\fA)_{1_{G^m}}$ leaves  $G_i$ fixed (as set). Moreover, by virtue of~\eqref{311222a}, we have $\aut(\fA)\le\aut(\Gamma)$, where $\Gamma=\Gamma_D(G,m)$. By Lemma~\ref{220522b}, this yields $\aut(\fA)_{1_{G^m}}\le\aut_m(G)$. Thus, 
$$
\aut(\fA)=G^m\aut(\fA)_{1_{G^m}}\le G^m\aut_m(G)=\hol_m(G).
$$
To prove the converse inclusion, it suffices to verify that $\aut_m(G)\le \aut(\fA)$. Every $f\in\aut_m(G)$ is a normalized automorphism of~$\cT_m$, such that $(D_mx)^f=D_mx^f$ for all $x\in G^m$. Thus, $f\in\aut(\fA)$ by the second part of Lemma~\ref{140514a}.
\eprf

\subsection{Algebraic isomorphism}
Let $\fA\ge \fA_m(G)$ and $\fA'\ge \fA_m(G')$. Then $G^{}_i$ and $G'_i$ are $\fA$- and $\fA'$-groups for all $i\in \ov M$. An algebraic isomorphism $\varphi:\fA\to \fA'$ is said to be \emph{genuine} if $\varphi(G_i)=G'_i$ for all~$i$; in particular,
\qtnl{070522b1}
\varphi(\cT_m(G))=\cT_m(G')\qaq \varphi(D_m(G))=D_m(G').
\eqtn
Clearly, $\varphi$ induces by restriction a genuine algebraic isomorphism from $\fA_m(G)$ to~$\fA_m(G')$. The following statement shows that every algebraic isomorphism of~$\fA_m(G)$ is genuine.

\prpstnl{280122d1}
Let $m\ge 3$, $\fA'$ an S-ring over a group $H$, and $\varphi:\fA_m(G)\to\fA'$ an algebraic isomorphism.  Then $H=(G')^m$ for a certain group $G'$. Moreover, $\fA'=\fA_m(G')$ and $\varphi$ is genuine. 
\eprpstn
\prf 
Note that  $G^m$ is the direct product of any $m$ of the $\fA_m(G)$-subgroups $G_i$, $i\in \ov M$.  Hence, $H$ is the direct product of any $m$ of $\fA'$-subgroups $H_i:=\varphi(G_i)$, $i\in\ov M$. By Lemma~\ref{220522a}, there is a group $G'$  such that $H=(G')^m$, $H_0=D_m(G')$ and $H^{}_i=G'_i$ for  all $i\in M$. It follows that $\fA'=\fA_m(G')$.
\eprf

We complete the subsection by collecting some properties of genuine isomorphisms, that will be used in Section~\ref{220522i}.

\lmml{280122d}
Let $m\ge 1$ and $\varphi$ a genuine algebraic isomorphism from $\fA_m$ to another S-ring. Then for all $X\in\cS_m$, $\sigma\in\sym(M)$, and $K\subseteq M$, we have
\nmrt
\tm{1} $\varphi(X^\sigma)=\varphi(X)^\sigma$,
\tm{2}  $\varphi_K(\pr_K(X))=\pr_K(\varphi(X))$, where $\varphi_K=\varphi_{G^m/G_{K'}}$,
\tm{3}  $n_K(X)=n_K(\varphi(X))$.
\enmrt
\elmm
\prf Denote by $\psi_\sigma$ the algebraic automorphism of $\fA_m$, induced by a permutation~$\sigma\in\sym(M)$. Then  the composition $\psi'=\varphi\psi_\sigma\varphi^{-1}$ is an algebraic automorphism of~$\fA_m(G')=\im(\varphi)$. By formulas~\eqref{070522b1} and Lemma~\ref{100522b},  we conclude that $\psi'$ is the algebraic automorphism of $\fA_m(G')$, induced by~$\sigma$. Thus, 
$$
\varphi(X^\sigma)=\varphi(\psi_\sigma(X))=\psi'(\varphi(X))=\varphi(X)^\sigma,
$$
which proves statement~(1).  Statement~(2) follows from formula~\eqref{180522a} for the natural epimorphism $\pi:G^m\to G^m/G_{K'}$. Finally, $n_K(X)= n(X,G_{K'})$, see the proof of Lemma~\ref{160422w}, and statement~(3) follows from formula~\eqref{180522h}.\eprf

\subsection{Isomorphisms}
In this subsection, we prove that the S-ring~$\fA_m(G)$ determines the group $G$ up to isomorphism.

\lmml{220522q}
Let $m\ge 1$. Assume that $G^m=H^m$. Then $G\cong H$. Moreover, if $\fA_m(G)=\fA_m(H)$,  then  there is a normalized isomorphism of~$\fA_m(G)$ that takes~ $G_i$ to~$H_i$ for all $i\in\ov M$.
\elmm
\prf
Let  $\{L_1,\ldots,L_k\}$ be a full set of pairwise nonisomorphic indecomposable (into direct product) subgroups of $G$ and hence of $G^m$. Since $H^m=G^m$, the Krull-Schmidt theorem implies that there are exactly $k$  pairwise nonisomorphic indecomposable subgroups $M_1,\ldots,M_k$ of $H^m$ such that $L_i\cong M_i$ for all~$i$. 

Note that if $n_i$  (respectively, $n'_i$) is the multiplicity of~$L_i$ (respectively, $M_i$) in a decomposition of $G_1$ (respectively, $H_1$) into indecomposable direct product, then  the multiplicity of~$L_i$ (respectively, $M_i$) in a decomposition of $G^m$ into indecomposable direct product is equal to $mn_i$ (respectively, $mn'_i$). Therefore, $mn^{}_i=mn'_i$ and hence $n^{}_i=n'_i$ for all~$i$. This shows that $G_i\cong H_i$ for $i=1$ and hence for all~$i$; in particular, $G\cong H$.  

The isomorphisms $G_i\to H_i$ induce an automorphism $f\in\aut(G^m)$  such that $H_i=(G_i)^f$ for all $i\in\ov M$. Since $\fA_m(G)$ is generated by both the $\und{G_i}$ and $\und{H_i}$, we have $f\in\iso(\fA_m(G))$.
\eprf

\crllrl{220522n}
Let $m\ge 3$. For any algebraic isomorphism  $\varphi:\fA_m(G)\to\fA_m(G')$ there exists $f\in\iso(\fA_m(G'))$ such that the composition $\varphi\varphi_f$  is a genuine algebraic isomorphism. 
\ecrllr
\prf 
By Proposition~\ref{280122d1}, there is a group $\wt G'$  such that $\fA_m(G')=\im(\varphi)=
\fA_m(\wt G')$ and the algebraic isomorphism $\varphi:\fA_m(G)\to\fA_m(\wt G')$ is genuine. 
 By Lemma~\ref{220522q}, there exists $f\in\iso(\fA_m(G'))$ such that $G'_i=(\wt G'_i)^f$ for all $i\in\ov M$. Thus,  the algebraic isomorphism $\varphi\varphi_f$ is genuine.
 \eprf

\crllrl{160522a}
Let $m\ge 3$. Then for any groups $G$ and $G'$, either all or none of algebraic isomorphisms from $\fA_m(G)$ to $\fA_m(G')$ is induced by an isomorphism.
\ecrllr
\prf Assume that an algebraic isomorphism $\varphi:\fA_m(G)\to\fA_m(G')$ is induced by an isomorphism. If $\varphi':\fA_m(G)\to\fA_m(G')$  is another algebraic isomorphism, then $\varphi^{-1}\varphi'$ is induced by an isomorphism by Lemma~\ref{220522q}. It follows that $\varphi'$ is also induced by an isomorphism. 
\eprf

{\bf Proof of Theorem \ref{220522t}.} It suffices to verify implications $(2)\Rightarrow(1)$ and $(3)\Rightarrow(1)$. Assume that  the S-rings $\fA_m(G)$ and $\fA_m(G')$ are isomorphic. Then by Corollaries~\ref{220522n} and \ref{160522a}, there is an isomorphism $f:\fA_m(G)\to \fA_m(G')$ inducing a genuine algebraic isomorphism. Without loss of generality, we may assume that~$f$ is normalized. Then $X_m(G)^f=X_m(G')$.  It follows that $f$ is a graph  isomorphism from $\Gamma_D(G,m)$ to $\Gamma_D(G',m)$. Thus, $G\cong G'$ by formula~\eqref{220522c}.  This proves implication $(2)\Rightarrow(1)$.

Assume that the groups $\hol_m(G)$ and $\hol_m(G')$ are permutation isomorphic. The permutation isomorphism is also an isomorphism of S-rings $\cyc_m(G)$ and $ \cyc_m(G')$, associated with that groups. In its turn, this isomorphism induces an S-ring isomorphism $f$ from $\fA_m(G)\le \cyc_m(G)$ to an S-ring $\fA'\le\cyc_m(G')$. By Proposition~\ref{280122d1} for the algebraic isomorphism induced by~$f$, we have $\fA'=\fA_m(\wt G)$ for some group $\wt G$ such that $(G')^m=(\wt G)^m$. By Lemma~\ref{220522q}, $\wt G\cong G'$ and hence $\fA'=\fA_m(G')$. It follows that the S-rings $\fA_m(G)$ and $\fA_m(G')$ are isomorphic and the implication $(3)\Rightarrow(1)$ is a consequence of the implication $(2)\Rightarrow(1)$.
\eprff\medskip

The arguments of this section can be used to strengthen the main result in~\cite{Rode2019}, stating that every group~$G$ is determined by the S-ring $\cyc_3(\Inn(G))^{\sym(3)}$ consisting of all elements of~$\cyc_3(\Inn(G))$ leaving fixed with respect to isomorphisms induced by the elements of $\sym(3)$. In fact, the latter S-ring can be replaced by a smaller S-ring, namely, $\fA_3(G)^{\sym(3)}$.

\section{Relationship between  $m$-ary coherent configurations  and S-rings}\label{220522i}

In this section, we establish two reductions between the $m$-ary coherent configurations $\WL_m(G)$ and S-rings $\fA_m(G)$, that will be used in the proof of the main results. 

\thrml{010422a}
Let  $m\ge 1$ and $\fA=\fA(G)$ the linear space defined by equality~\eqref{300322a} for $\cS=\pr_m(\WL_{3m}(G))$. Then
\nmrt
\tm{1} $\fA(G)$ is an S-ring over~$G^m$,
\tm{2} $\fA(G)\ge \fA_m(G)$,
\tm{3} if $G'$ is a group such that $\WL_{3m}(G)$ and $\WL_{3m}(G')$ are genuine algebraically isomorphic, then so are $\fA(G)$ and $\fA(G')$.
\enmrt
\ethrm
\prf By statements~(1) and (2) of Lemma~\ref{140422b}, the partition $\WL_{3m}(G)\ge \fX_{3m}(G)$ contains the class $\{1_{G^{3m}}\}$ and is closed with respect to taking the inverse. Since the projection $\pr_m:G^{3m}\to G^m$ is a group homomorphism,  $\cS$ contains the class $\{1_{G^m}\}$ and is closed with respect to taking the inverse. Thus $\cS$ satisfies the conditions~(S1) and  (S2). 

To verify the condition (S3), let $X,Y,Z\in\cS$. We need to find an integer $c=c_{X,Y}^Z$ such that every $z\in Z$ has exactly $c$ representations $xy=z$ with $x\in X$ and $y\in Y$. To this end, put $A=X\times Y\times Z$. Then $A\subseteq G^{3m}$ and
$$
\pr_{1,\ldots,m}(A)=X,\qquad
\pr_{m+1,\ldots,2m}(A)=Y,\qquad 
\pr_{2m+1,\ldots,3m}(A)=Z.
$$
Since $\fX=\WL_{3m}(G)$ is a $3m$-ary coherent configuration, formula~\inp{\eqref{160523g}} shows that  the full preimages $X'$, $Y'$, and $Z'$ of the classes $X$, $Y$, and $Z$ with respect to the corresponding projections 
belong to $\fX^\cup$. Consequently,  $A=X'\cap Y'\cap Z'$ belongs to $\fX^\cup$.  Furthermore, the right-hand side of~\eqref{030422a} implies that $\fX^\cup$ contains also the set 
\qtnl{160522u}
A'=\{(x,y,z)\in G^{3m}:\ x_iy_i=z_i,\quad i=1,\ldots,m\}.
\eqtn
Therefore, $A\cap A'\in \fX^\cup$ and $A\cap A'=A_1\cup\ldots\cup A_r$ for some $A_1,\ldots,A_r\in \fX$ and $r\ge 0$.  Since $r=0$ if and only if $XY\cap Z=\varnothing$,  we may assume that $r\ge 1$. 

Let $z\in Z$. Since $r\ge 1$, there is at least one pair  $(x,y)\in X\times Y$ such that $z=xy$. For any such a pair,  $(x,y,z)\in A'$, and hence $(x,y,z)\in A_i$ for exactly one $i\in\{1,\ldots, r\}$. By Lemma~\ref{090122d}, the number of all $(x,y)\in X\times Y$ for which $(x,y,z)\in A_i$ is equal to the number $n_K(A_i)$ with $K=\{2m+1,\ldots,3m\}$. Thus, the number
\qtnl{110522f}
c_{X,Y}^Z=c=|\{(x,y)\in X\times Y:\ xy=z\}|=\sum_{i=1}^rn_K(A_i)
\eqtn
does not depend on $z\in Z$, as required. This completes the proof of statement~(1).

Furthermore, by Corollary~\ref{180422b} and Lemma~\ref{140422b}(3), we have
$$
\WL_{3m}(G)\ge \fX_{3m}(G)\ge \cS(\cT_{3m})\qaq D_{3m}\in \fX_{3m}(G)^\cup\le \WL_{3m}(G)^\cup.
$$
It follows that $\cS\ge  \pr_m(\cS(\cT_{3m}))=\cT_m$ and $D_m=\pr_m(D_{3m})\in \cS^\cup$. Thus the S-ring~ $\fA(G)$ contains the extension of $\cT_m$ by $D_m$, i.e., $\fA_m(G)$. This proves statement~(2).

To prove statement~(3), let $\varphi:\WL_{3m}(G)\to\WL_{3m}(G')$ be a genuine algebraic isomorphism. By the second part of Lemma~\ref{090122b}, it induces a uniquely determined bijection $\psi:\cS(\fA(G))\to \cS(\fA(G'))$ such that for all $\wt A\in\WL_{3m}(G)$,
\qtnl{160522t}
 \psi(\pr_m(\wt A))=\pr_m(\varphi(\wt A)).
 \eqtn
Now let  $X,Y,Z\in\cS$ and $A=X\times Y\times Z$. Then $\varphi(A)=\psi(X)\times\psi(Y)\times\psi(Z)$. Since the algebraic isomorphism $\varphi$ is genuine,  formula~\eqref{100522j} yields 
$$
\mu(\varphi(A'))=\mu(A')\supseteq\{(i,i+m,i+2m):\ i\in M\},
$$ 
where $A'$ is defined by formula~\eqref{160522u}.  Thus, as above, $A\cap A'=A_1\cup\ldots\cup A_r$ and using Lemma~\ref{280122d}(3), we obtain
$$
c_{\psi(X),\psi(Y)}^{\psi(Z)}=\sum_{i=1}^rn_K(\inp{\varphi}(A_i))=\sum_{i=1}^rn_K(A_i)=c_{X,Y}^Z.
$$
It follows that $\psi$ induces an algebraic isomorphism from $\fA(G)$ to $\fA(G')$. It is genuine, because so is~$\varphi$. 
\eprf

\thrml{030422s}  
Let $m\ge 1$ and $\fX(G)=\pr_m(\cS_{m+1})$. Then
\nmrt
\tm{1} $\fX(G)$  is an $m$-ary coherent configuration,
\tm{2} $\fX(G)\ge \WL_m(G)$ if $m\ge 2$,
\tm{3} if $G'$ is a group such that $\fA_{m+1}(G)$ and $\fA_{m+1}(G')$ are genuine algebraic isomorphic, then so are $\fX(G)$ and $\fX(G')$.
\enmrt
\ethrm

\prf  By Lemma~\ref{160422w}, the partition $\cS_{m+1}$ is a regular $(m+1)$-ary rainbow. Therefore the partition $\fX(G)$ is an $m$-ary rainbow by Lemma~\ref{090122b}. Let us verify that $\fX=\fX(G)$ satisfies condition~(C3). 

Let $X_0=:X,X_1\ldots,X_m\in \fX$. By formula~\eqref{160523g}, the full preimage $Y=X\times G$ of~$X$ with respect to~$\pr_m$ belongs to $(\cS_{m+1})^\cup$. It follows that so does
$$
Y_i=\{y\in Y:\ \pr_m(y^{\sigma_i})\in X_i\},\qquad i\in M,
$$
where $\sigma_i=(i,m+1)$ is a transposition of $\sym(m+1)$. Consequently the set $A=Y_1\cap\ldots\cap Y_m$ belong to~$\cS_{m+1}$. Note that if $(x,g)\in Y$ and $\wh x=(x_1,\ldots,x_m,g)$, then $x_{i\leftarrow g}=\pr_m(\wh x^{\sigma_i})$ and
\qtnl{120122a}
A=\{(x,g)\in Y:\ x_{i\leftarrow g}\in X_i,\ i\in M\}.
\eqtn
Let $A$ be the disjoint union of $A_1,\ldots,A_r\in \cS_{m+1}$, where $r\ge 0$. Since the $(m+1)$-ary rainbow~$\cS_{m+1}$ is regular,  the number $n_{\{m+1\}}(A_i)$ of all $g\in G$ such that $(x,g)\in A_i$ does not depend on $x\in\pr_m(A_i)$,   $i=1,\ldots,r$. Thus by formula~\eqref{120122a},  the number
$$
n_{X_1,\ldots,X_m}^X=|\{g\in G:\  x_{j\leftarrow g}\in X_j,\ j=1,\ldots,m\} |=
$$
$$
|\{g\in G:\ (x,g)\in A\}| =\sum_{i=1}^r|\{g\in G:\ (x,g)\in A_i\}| =\sum_{i=1}^r n_{m+1}(A_i)
$$
also does not depend on $x\in\pr_m(A_i)=X$.  Consequently the $m$-ary rainbow $\cS_m$ satisfies condition (C3) and hence is an $m$-ary coherent configuration. This proves statement~(1).

By  Corollary \ref{180422b}, we have $\cS_{m+1}\ge \fX_{m+1}(G)$. After taking the projection to~$G^m$ and using Lemma~\ref{140422b}(3), we obtain  
$$
\fX(G)=\pr_m(\cS_{m+1})\ge \pr_m(\fX_{m+1}(G))\inp{\ge}\fX_m(G).
$$
By statement~(1) and  the monotonicity  of the operator~$\WL_m$, we obtain
$$
\fX(G)=\WL_m(\fX(G))\ge \WL_m(\fX_m(G))=\WL_m(G),
$$
which proves statement~(2).

To prove statement~(3), let  $\varphi:\fA_{m+1}(G)\to\fA_{m+1}(G')$ be a  genuine algebraic isomorphism.  The classes of $\fX(G)$ and $\fX(G')$ are the basic sets of $\fA_{m+1}(G)_{G^m}$ and $\fA_{m+1}(G')_{(G')^m}$, respectively. Therefore the isomorphism~$\varphi$ induces a bijection $\psi:\fX(G)\to\fX(G')$ such that $\varphi(X\times G)=\psi(X)\times G'$ for all $X\in\fX(G)$. By statements~(1) and~(2) of Lemma~\ref{280122d}, we have
$$
\varphi(Y_i)=\{y\in \varphi(Y):\ \pr_m(y^{\sigma_i})\in\psi(X_i)\},\qquad i\in M,
$$
where $Y_i$ and $X_i$  are as above. It follows that $n_{\psi(X_1),\ldots,\psi(X_m)}^{\psi(X)}=\sum_{j=1}^r n_{m+1}(\psi(A_j))$. By statement~(3) of Lemma~\ref{280122d}, this yields
$$
n_{\psi(X_1),\ldots,\psi(X_m)}^{\psi(X)}=\sum_{i=1}^r n_{m+1}(\psi(A_i))=\sum_{i=1}^r n_{m+1}(A_i)=n_{X_1,\ldots,X_m}^X.
$$
Consequently, $\psi$ is an algebraic isomorphism from $\fX(G)$ to $\fX(G')$, which is genuine because so is~$\varphi$.
\eprf

\section{Proof of the main results }\label{190122c}

\subsection*{Proof of Theorem~\ref{271221a0}}  Let us prove inclusions~\eqref{300322g}. It is easily seen that
$$
\pr_m(\cT_{m+k})=\cT_m\qaq \pr_m(D_{m+k})=D_m.
$$
Hence the S-ring $\fA=(\fA_{m+k})_{G^m}$ contains the extension of the S-ring $\cT_m\le \fA$ by the set $D_m\in \cS(\fA)^\cup$. Thus, $\fA\ge \fA_m$, which proves the first inclusion in~\eqref{300322g}. Next, $\aut(\fA_{m+k})=\hol_{m+k}(G)$ by Corollary~\ref{220522w}. Therefore, $\fA_{m+k}\le \cyc_{m+k}(G)$ and hence
$$
(\fA_{m+k})_{G^m}\le (\cyc_{m+k}(G))_{G^m}=\cyc_m(G).
$$
This proves the second inclusion in~\eqref{300322g}.

Assume that $k\ge \max\{2,d\}$. Then $k+m\ge 3$.  Let  $x\in G^{m+k}$ be an arbitrary element such that $\grp{x_1,\ldots,x_d}=G$. By Corollary~\ref{290422f}, for each $y\in [x]$, there is a group isomorphism
$$
\sigma_y:\grp{x_1,\ldots,x_{m+k}}\to\grp{y_1,\ldots,y_{m+k}}, \ x_i\mapsto y_i\ (i=1,\ldots m+k).
$$
Then $\sigma_y\in \aut(G)$ and $y=x^{\sigma_y}$.  Since also $X=[x]$ is invariant with respect to~$\aut_{k+m}(G)$, this implies that $X$ is an orbit of $\aut_{k+m}(G)$ and hence $X'=\pr_{k+1,\ldots,m+k}(X)$ is an orbit of $\aut_m(G)$. When  the tuple $\pr_{k+1,\ldots,m+k}(x)$ runs over the group~$G^m$, the projection $X'$ runs over the basic sets of~$\fA$. Therefore $\cS(\fA)$ consists of the orbits of  $\aut_m(G)$, and $\fA=\cyc(\aut(G),G^m)=\cyc_m(G)$.

\subsection*{Proof of Theorem~\ref{230122a}}  By Theorem~\ref{010422a}(1), the partition  $\pr_m(\WL_{3m}(G))$ is equal to the partition $\cS(\fA)=:\cS$ with $\fA=\fA(G)$. Moreover, by Theorem~\ref{010422a}(2), $\fA\ge \fA_m$, in particular, $\cS\ge \cS(\fA_m)$. Thus,
$$
\pr_m(\WL_{3m}(G))=\cS\ge \cS(\fA_m(G)),
$$
which  proves the first inclusion in Theorem~\ref{230122a}. Similarly, by Theorem~\ref{030422s}(1), the partition  $\pr_m(\cS_{m+1})$ is an $m$-ary coherent configuration. Moreover, by Theorem~\ref{030422s}(2), $\fX\ge \WL_m(G)$, which  proves the  second inclusion in Theorem~\ref{230122a}. 

Let us prove that if $m\ge \wldim(G)+1$,  then every algebraic isomorphism from the S-ring  $\fA_m(G)$ is induced by an isomorphism; this proves the existence of~$s(G)$ and  the second inequality in~\eqref{090123a}.

Let $m=\wldim(G)$, $m\ge 3$, and let $G'$ be a group such that the S-rings $\fA_{m+1}(G)$ and $\fA_{m+1}(G')$ are algebraically isomorphic. By Corollary~\ref{220522n}, we may assume that they are genuine algebraically isomorphic. By statements~(2) and~(3) of Theorem~\ref{030422s}, this 
implies that the $m$-ary coherent configurations $\fX(G)\ge \WL_m(G)$ and $\fX(G')\ge \WL_m(G')$ are genuine algebraically isomorphic.  By Lemma~\ref{130522v} for $\fX=\fX(G)$ and $\fY=\WL_m(G)$, so are the $m$-ary coherent configurations $\WL_m(G)$ and $\WL_m(G')$. Since  $m=\wldim(G)$, we conclude that $G\cong G'$, and the S-rings $\fA_{m+1}(G)$ and $\fA_{m+1}(G')$ are isomorphic  by Theorem~\ref{220522t}. In particular, a genuine algebraic isomorphism between them is induced by an isomorphism. By Corollary~\ref{160522a}, this shows that every algebraic isomorphism from the S-ring  $\fA_{m+1}(G)$ is induced by an isomorphism. This proves the existence of the number~$s(G)$ and  the second inequality in~\eqref{090123a}.

Let $m=s(G)$, $m\ge 3$, and let $G'$ be a group  $\WL_{3m}$-equivalent to~$G$.  This means that the $m$-ary coherent configurations $\WL_{3m}(G)$ and $\WL_{3m}(G')$ are genuine algebraically isomorphic. By statements~(2) and~(3) of Theorem~\ref{010422a}, this implies that so are the S-rings $\fA(G)\ge \fA_m(G)$ and $\fA(G')\ge \fA_m(G')$. It follows that the S-rings $\fA_m(G)$ and $\fA_m(G')$ are algebraic isomorphic. Since  $m=s(G)$, they are isomorphic, and hence $G\cong G'$ by Theorem~\ref{220522t}. Thus, $\wldim(G)\le 3s(G)$.

\section{Proof of Theorem~\ref{220122e}}\label{041122a}

The  statement of Theorem~\ref{220122e} is an immediate consequence of a more general lemma below. In  the proof, under the individualization $G_x$ of a colored group $G$ by an  element~$x\in G$, we mean that~$G$ is colored so that the colors of the elements~$\ne x$ are the same as before, whereas the color of~$x$ is different from the colors of the other elements.

\lmml{220122e1}
The following problems for colored groups $G$ and $G'$ are polynomial-time equivalent:
\nmrt
\tm{a} test whether $\iso(G,G')\ne\varnothing$, and if so find an element of $\iso(G,G')$,
\tm{b} find the set $\iso(G,G')$,
\tm{c} find the group $\aut(G)$,
\tm{d} find the S-ring $\cyc_1(G)$.
\enmrt
\elmm
\prf (a) reduces to (b) trivially. To describe the reduction of~(b) to~(c), denote by  $c_{G^{}}$ and $c_{G'}$ the colorings of~$G$ and $G'$, respectively. We define a coloring~$c$  of the direct product $G\times G'$ so that for any $g\in G$ and $g'\in G'$,
\qtnl{300122a}
c(g,g')=\css
c_{G^{}}(g)   &\text{if $g\ne 1$ and $g'=1$,}\\
c_{G'}(g')   &\text{if $g=1$ and $g'\ne 1$,}\\
\varepsilon &\text{otherwise,}\\
\ecss
\eqtn
where $\varepsilon$ is an arbitrary symbol not in $\im(c_{G^{}})\cup\im(c_{G'})$. It is easily seen that the automorphism group $\aut(G\times G')$ of the colored group $G\times G'$ leaves the set~$G\cup G'$ fixed; denote by $\aut_0(G\times G')$ the subgroup of $\aut(G\times G')$ leaving both $G$ and $G'$ fixed (as sets). It remains to note that $G\cong G'$ if and only if 
$$
[\aut(G\times G'):\aut_0(G\times G')]=2,
$$
and if so, then $\iso(G,G')=Kf$, where $K$ is the restriction of $\aut_0(G\times G')$ to~$G$ and $f:G\to G'$ is a bijection induced by an (arbitrarily chosen) permutation belonging to the nontrivial coset of $\aut_0(G\times G')$ in $ \aut(G\times G')$.

Let us describe a reduction of~(c) to~(a). Assume that we are given an algorithm solving~(a). Then for every element~$x$ of a colored group $G$, one can efficiently find
the set 
$$
S_x(G)=\{f_{x,y}\in \iso(G_x,G_y):\ y\in G,\ G_x\cong G_y\},
$$
where for each~$y$ the isomorphism $f_{x,y}$ is chosen arbitrarily. Now if the coloring of $G$ is
discrete, then, of course, $\aut(G)$ is trivial. Next, if  $\{x\}$ is not a color class of $G$ for at least one  $x\in G$, then $\aut(G)=\grp{S_x(G),\aut(G_x)}$.
Thus finding $\aut(G)$ is efficiently reduced to finding the groups $\aut(G_x)$. Moreover, the number of singleton color classes of~$G_x$ is at most $|G|$.  Consequently the group $G$ can be constructed in at most $|G|$ reductions.  

(d)  reduces to ~(c) trivially. Let us describe a reduction of~(a) to~(d). We need an auxiliary statement.\medskip

{\bf Claim.} {\it Let $G$ and $G'$ be colored groups. Then given $x\in G$ and $x'\in G'$,  one can efficiently test by using (d) as oracle whether or not $\iso(G^{}_{x^{}},G'_{x'})\ne\varnothing$.}

\prf
Define a coloring of  $K=G^{}_{x^{}}\times G'_{x'}$ by formula~\eqref{300122a} (note that the pairs $(x,1)$ and $(1,x')$ are colored in the same color). Then  $\iso(G^{}_{x^{}},G'_{x'})\ne\varnothing$ if and only if $x$ and $x'$ lie in the same orbit of $\aut(K)$ if and only if  $(x,1)$ and $(1,x')$ lie in the same  class of the partition $\cS(\fA)$ with $\fA=\cyc_1(\aut(K))$.
\eprf

Now let $G$ and $G'$ be colored groups. If the coloring of $G$ is discrete, then~(a) is solved in an obvious way. Otherwise, let $x\in G$.  By the claim, one can efficiently test whether or not $\iso(G^{}_{x^{}},G'_{x'})\ne\varnothing$ for at least one~$x'$. If not, then, of course, $\iso(G,G')=\varnothing$. Otherwise, the problem~(a) for $G$ and~$G'$ is efficiently reduced to that for~$G^{}_{x^{}}$ and~$G'_{x'}$. Moreover, the number of singleton color classes of~$G_x$ is strictly less than that of~$G$.  Consequently after at most $|G|$ reductions, we arrive at the problem~(a) in which one of the input groups has discrete coloring.
\eprf

\section*{Acknowledgment}
The work of the first author is supported by Natural Science Foundation of China (No. 11971189, No. 12161035)


\begin{thebibliography}{10}
	
\bibitem{Babai2019}
L.~Babai, \emph{{Group, Graphs, Algorithms: the Graph Isomorphism Problem}}, 	Proceedings of the International Congress of Mathematicians (ICM 2018), 	vol.~3, WORLD SCIENTIFIC (2019), see also L.~Babai, {\em Graph Isomorphism in 		Quasipolynomial Time} (2016), arXiv:1512.03547v2 [cs.DS]), pp.~3319--3336.
	
\bibitem{Baginski2021}
C.~Bagi{\'{n}}ski and P.~Grzeszczuk, \emph{{On the generic family of Cayley graphs of a finite group}}, J. Combin. Theory. Ser. A, \textbf{184},	105495 (2021).
	
\bibitem{Bailey2021}
R.~A. Bailey, P.~J. Cameron, C.~Praeger, and C.~Schneider, \emph{{The geometry 	of diagonal groups}}, Trans. Amer. Math. Soc., {\bf 375}, no. 8, 5259--5311 (2022).
	
\bibitem{Brachter2020}
J.~Brachter and P.~Schweitzer, \emph{{On the Weisfeiler-Leman Dimension of 	Finite Groups}}, Proc. 35th Annual ACM/IEEE Symposium on Logic in Computer Science (New York, NY, USA), no.~1, ACM, 2020, pp.~287--300.
	
\bibitem{Brachter2021}
J.~Bra{ch}ter and P.~Schweitzer, \emph{{A Systematic Study of Isomorphism 	Invariants of Finite Groups via the Weisfeiler-Leman Dimension}}, in: \emph{30th Annual European Symposium on Algorithms (ESA 2022)}, Article No.~27 (2022), pp. 27:1--27:14.

\bibitem{CP2019}
G.~Chen and I.~Ponomarenko, \emph{Coherent configurations}, Central China 	Normal University Press (2019), the updated version is available at \url{http://www.pdmi.ras.ru/~inp/ccNOTES.pdf}.
	
\bibitem{EvdP1999c}
S.~Evdokimov and I.~Ponomarenko, {\em On highly closed cellular algebras and highly closed isomorphisms}, Electronic J. Combin., {\bf 6}, \#R18 (1999).
	
\bibitem{Evdokimov2013}
S.~Evdokimov and I.~Ponomarenko, \emph{Schur rings over a product of {G}alois rings}, Beitr. Algebra Geom., \textbf{55}, no.~1, 105--138  (2014).
	
\bibitem{Evdokimov2015}
S.~Evdokim{o}v and I.~Ponomarenko, \emph{{On the separability problem for 	circulant S-rings}}, St. Petersburg Math. J., \textbf{28}, no.~1, 21--35  (2017).
	
\bibitem{Formanek1991}
E.~Formanek and D.~Sibley, \emph{{The Group Determinant Determines the Group}}, Proc. AMS, \textbf{112}, no.~3, 649--656  (1991).
	
\bibitem{Grechkoseeva2021}
M.~Grechkoseeva, V.~D. Mazurov, W.~Shi, A.~Vasil'ev, and N.~Yang, \emph{{Finite groups isospectral to simple groups}}, Commun. Math. Stat. (2022), doi: 10.1007/s40304-022-00288-5.
	
\bibitem{Heinze2009}
A.~Heinze and M.~Klin, \emph{{Loops, Latin Squares and Strongly Regular Graphs: An Algorithmic Approach via Algebraic Combinatorics}}, in: Algorithmic Algebraic 	Combinatorics and Gr{\"{o}}bner Bases, Springer Berlin Heidelberg, Berlin, 	Heidelberg (2009), pp.~3--65.
	
\bibitem{AndresHelfgott2017}
H.~Helfgott, J.~Bajpai, and D.~Dona, \emph{{Graph isomorphisms in quasi-polynomial time}}, {\tt arXiv:1710.04574}, 1--67 (2017).
	
\bibitem{Mat1979}
R.~Mathon, \emph{A note on the graph isomorphism counting problem}, Inform. Process. Lett. \textbf{8}, 131--132 (1979).
	
\bibitem{Ponomarenko2022a}
I.~Pono{m}arenko, \emph{{On the WL-dimension of circulant graphs of prime power order}}, {\tt arXiv:2206.15028}, 1--23 (2022).
	
\bibitem{Rode2019}
E.~L. Rode, \emph{{On a generalized centralizer ring of a finite group which determines the group}}, Algebra Colloq., \textbf{26}, no.~1, 31--50  (2019).
	
\bibitem{Roitman1981}
M.~Roitman, \emph{{A complete set of invariants for finite groups and other results}}, Advances Math. \textbf{41}, no.~3, 301--311  (1981).
	
\bibitem{Wielandt1964}
H.~Wielandt, \emph{Finite permutation groups}, Academic Press, New York and London (1964).
	
\end{thebibliography}

\end{document}